 \documentclass[11pt]{amsart}
\usepackage{amscd,amssymb,euscript,amsthm}
\theoremstyle{plain}
\newtheorem{thm}{Theorem}[section] 
\newtheorem{cor}[thm]{Corollary}
\newtheorem{prop}[thm]{Proposition}
\newtheorem{lem}[thm]{Lemma}
\theoremstyle{definition}
\newtheorem{defn}[thm]{Definition}
\theoremstyle{remark}
\newtheorem{rem}[thm]{Remark}

\numberwithin{equation}{section}%assigns eqns section numbers
\newcommand{\ball}{\operatorname{ball}}

\def\<{\left<}
\def\>{\right>}
\def\cstar{$C^*$-algebra}
\newcommand{\comment}[1]{}
\begin{document}
\title[Asymptotic lift]{Asymptotic lifts of  positive linear maps}
\author{William Arveson and Erling St\o rmer}
%
%\thanks{*supported by 
%NSF grant DMS-0100487} 
%
\address{Department of Mathematics,
University of California, Berkeley, CA 94720
Department of Mathematics, Univesity of Oslo, 
0316 Oslo, Norway}
\email{arveson@math.berkeley.edu, erlings@math.uio.no}
%
%\subjclass{46L55, 46L09}
%\date{18 June, 2006}
%
\begin{abstract} We show that the notion of asymptotic lift generalizes 
naturally to normal positive maps $\phi: M\to M$ acting 
on von Neumann algebras $M$.  
We focus on cases in which the domain of the asymptotic lift can be embedded 
as an operator subsystem $M_\infty\subseteq M$, and characterize 
when $M_\infty$ is a Jordan subalgebra of $M$ in terms of the 
asymptotic multiplicative  properties 
of $\phi$.  
\end{abstract}
\maketitle

\section{Introduction}\label{in}

Let $\phi:M\to M$ be a normal unit preserving positive 
linear map acting on a dual operator system $M$; we refer to such 
a pair $(M,\phi)$ as a {\em UP map}.  While we are primarily interested in 
UP maps that act on von Neumann algebras $M$, it is useful to 
broaden the context as above.  The powers of $\phi$ form an irreversible dynamical 
semigroup $\{\phi^n: n=0,1,2,\dots\}$ acting on $M$.  In this paper we generalize work begun in 
\cite{arvStableI} and \cite{arvLift}, together with complementary results 
in \cite{storPos}, to further develop the asymptotic theory of 
such semigroups.  

One may view UP maps as the objects of a category, 
in which a homomorphism from $\phi_1:M_1\to M_1$ to $\phi_2:M_2\to M_2$ is a 
UP map $E: M_1\to M_2$ such that $E\circ \phi_1=\phi_2\circ E$.  There is a natural 
notion of isomorphism in this category.  
Our first general result is that every normal unit-preserving positive linear map acting 
on a dual operator system has an 
asymptotic lift which is unique up to isomorphism.  This generalizes 
one of the main results of \cite{arvLift}, which dealt with the subcategory 
in which the objects are unital normal {\em completely} positive maps (UCP maps) 
on dual operator systems, with UCP maps as morphisms.  

We are primarily concerned with UP maps that act on von Neumann algebras $M$.  
In \cite{arvLift} it was shown that the asymptotic lift $(N,\alpha,E)$ 
of a UCP map $\phi:M\to M$ acting on a von Neumann algebra $M$ also acts on 
a von Neumann algebra $N$.  Moreover, it was shown that the $W^*$-dynamical 
system $(N,\alpha)$ can 
be identified with the tail flow of the minimal dilation 
of $\phi$ to a $*$-endomorphism of a larger von Neumann algebra in most cases -- namely  
those in which the 
dilation endomorphism has trivial kernel, which includes all 
UCP maps acting on factors $M$.  Since 
the minimal dilation of a UCP map on a von Neumann algebra can 
be constructed explicitly in principle, that 
provided a concrete identification of the asymptotic lift.  

It is significant that the asymptotic lift $(N,\alpha,E)$ of a UP map acting on a von 
Neumann algebra need not act on a von Neumann algebra $N$.  In this paper we show that 
in general, $N$ is order-isomorphic to a unique JW$^*$-algebra 
 -- namely a dual operator system that is closed under the Jordan product 
$x\circ y=(xy+yx)/2$ -- in such a way that $\alpha$ becomes a Jordan automorphism of $N$.  Thus, 
{\em the asymptotic behavior of a UP map on a von Neumann algebra 
is always associated with a JW$^{\, *}$-dynamical system $(N,\alpha)$}.  

Naturally, one would like to identify $(N,\alpha)$ more concretely in terms 
of the asymptotic properties of the semigroup $\phi,\phi^2,\phi^3,\dots$.  But 
since there is no dilation theory for semigroups of UP maps -- and perhaps an effective  
dilation theory is impossible, there seems to 
be no candidate to 
replace the ``tail flow" identification described above for UCP maps.  
Thus, this identification 
problem becomes a significant issue for UP maps.  

In many cases, $N$ can be embedded as an operator 
system $N\subseteq M$ in such a way that $\alpha=\phi\restriction_N$.  
While the operator system 
$N$ is order-isomorphic to a JW$^*$-algebra in general, it need not be 
closed under the ambient Jordan multiplication of $M$, and we 
address that issue in Section \ref{S:ec}.  We identify $N$ 
concretely in terms of the action of $\phi$ on $M$, and 
when $N$ is a Jordan subalgebra of $M$, we are able to go farther 
by identifying $N$  
with the {\em multiplicative core} of 
$\phi$ that was introduced in \cite{storPos}.  

The problem of finding a 
satisfactory concrete description of the JW$^*$-dynamical system 
$(N,\alpha)$ when $N$ does {\em not} embed in $M$ remains untouched.

\section{Asymptotic lifts of positive maps}\label{pr}

In this section we describe how the notion of asymptotic lift 
(of a completely positive map) can be generalized to normal 
unit-preserving maps that are 
merely positive.  We summarize 
the basic properties of asymptotic lifts,  indicating briefly 
how proofs of \cite{arvLift} should be modified.

\begin{rem}[Dual operator systems, norm, and order]\label{rem1}
Every von Neumann algebra is the dual of a unique Banach space, 
and hence it carries a natural weak$^*$ topology 
(aka the ultraweak or $\sigma$-weak topology).  
A {\em dual operator system} is a linear 
space of operators $M\subseteq \mathcal B(H)$ that contains 
the identity $\mathbf 1$, is self-adjoint $M^*=M$, and is 
closed in the weak$^*$-topology of $\mathcal B(H)$.  Such an 
$M$ also has a unique predual $M_*$, and its intrinsic $M_*$-topology 
coincides with the relative weak$^*$-topology of $\mathcal B(H)$.  
A map of dual operator systems 
$\phi:M\to N$ is called {\em normal} when it is weak$^*$-continuous.  
There is an intrinsic characterization of dual operator systems that 
we do not require.  

Let us recall the basic 
properties of unit-preserving normal positive linear maps 
$\phi: M\to N$ (UP maps) between dual operator systems.  
A UP map $\phi: M\to N$ defined on a 
\cstar\ $M$ must satisfy $\|\phi\|=1$.  That is most easily 
seen by making use of the Russo-Dye theorem \cite{rusDye} which implies that 
$\|\phi\|=\sup\|\phi(u)\|$, $u$ ranging over the unitary group of $M$, together with
the fact that 
for every unitary $u$, the 
restriction of $\phi$ to the commutative \cstar\ $C^*(u)$ is 
completely positive and therefore satisfies the strong Schwarz inequality 
$\phi(u)^*\phi(u)\leq \phi(u^*u)=\phi(\mathbf 1_M)=\mathbf 1_N$.  
More generally, if 
$M$ is merely an operator system then one has  
$1\leq \|\phi\|\leq 2$ in general, but the $C^*$-algebraic 
upper bound $\|\phi\|=1$ often fails.  

On the other hand, since the norm of a self-adjoint operator $x$ is the smallest  
$\alpha>0$ such that $-\alpha\mathbf 1\leq x\leq \alpha\mathbf 1$, the  restriction of $\phi$ to 
the real Banach space $M^{\rm sa}$ of self-adjoint elements of $M$ has norm $1$.  Conversely, 
if a linear map $\phi: M\to N$ carries self-adjoint elements to self-adjoint elements, maps  
$\mathbf 1_M$ to $\mathbf 1_N$, and satisfies $\|\phi\restriction_{M^{\rm sa}}\|=1$, 
then $\phi$ must also preserve positivity.  In particular,  {\em an order isomorphism 
in the category of UP maps is characterized as a $*$-preserving linear map $\phi: M\to N$ 
of operator systems that 
carries $\mathbf 1_M$ to $\mathbf 1_N$ 
and restricts to an isometry of $M^{\rm sa}$ onto $N^{\rm sa}$.} 
\end{rem}

\begin{defn}\label{inDef0}
A {\em reversible lift} of a UP map $\phi:M\to M$ is a triple $(N,\alpha,E)$ consisting 
of a UP automorphism $\alpha: N\to N$ of another dual operator system $N$ and a UP map 
$E: N\to M$ satisfying $E\circ\alpha=\phi\circ E$.  
\end{defn}

A reversible lift $(N,\alpha,E)$ of $\phi$ is said to be {\em nondegenerate} 
if 
\begin{equation}\label{inEq0}
 E(\alpha^{-n}(y))=0, \ n=0,1,2,\dots \ \implies y=0.  
\end{equation}
Significantly, if (\ref{inEq0}) fails, 
one can replace $(N,\alpha,E)$ 
with another
reversible lift $(\tilde N, \tilde\alpha,\tilde E)$ 
that is  {\em nondegenerate}, as in Remark 2.3 of \cite{arvLift}.

\begin{defn}\label{inDef1}
Let $\phi:M\to M$ be a UP map on a dual operator system.  
An {\em asymptotic lift} of $\phi$ 
is a reversible lift $(N,\alpha, E)$ of $\phi$ 
that satisfies nondegeneracy (\ref{inEq0}), 
together with 
\begin{equation}\label{inEqA}
\|\rho\circ  E\restriction_{N^{\rm sa}}\|=
\lim_{k\to\infty}\|\rho\circ  \phi^k\restriction_{M^{\rm sa}}\|,\qquad 
\rho\in M_*.    
\end{equation}
\end{defn}

\begin{rem}\label{inRem1}
We shall make use of the dual formulation of (\ref{inEqA}), and we record that now for later 
reference:  For every nondegenerate reversible lift $(N,\alpha,E)$ of $\phi$, (\ref{inEqA}) 
is equivalent to the following assertion:
\begin{equation}\label{inEqB}
E(\ball_r N^{\rm{sa}})=\bigcap_{n=1}^\infty \phi^n(\ball_r M^{\rm{sa}}), \qquad r>0, 
\end{equation}
where $\ball_r X$ denotes the closed ball of radius $r$ in a 
real or complex Banach space $X$.  
The proof of equivalence of (\ref{inEqA}) and (\ref{inEqB}) follows 
the lines of the proof of the corresponding result 
of \cite{arvLift}.  

\end{rem}

There are two fundamental results on asymptotic lifts of UP maps.  
The first concerns 
existence and uniqueness:  

\begin{thm}\label{inThm1}
Every UP map $\phi:M\to M$ of a dual operator system has an asymptotic lift.    
If $(N_1, \alpha_1, E_1)$ and 
$(N_2, \alpha_2, E_2)$ are two asymptotic lifts 
for $\phi$, then there is a 
unique UP-isomorphism of dual operator systems $\theta: N_1\to N_2$ such that 
$\theta\circ\alpha_1=\alpha_2\circ\theta$ and $E_2\circ\theta=E_1$.  
\end{thm}

As in the case of completely positive maps, 
the existence issue is settled by a direct construction based on {\em inverse sequences}, 
namely bounded bilateral sequences $(x_n)$ of elements of $M$ that satisfy  
$x_n=\phi(x_{n+1})$, $n\in\mathbb Z$.  The space of all inverse sequences 
is a dual operator system $N$, the bilateral shift $\alpha: (x_n)\mapsto (x_{n-1})$ 
is an automorphism of $N$, and the connecting map $E$ carries $(x_n)$ to $x_0$.  The 
proof that $(N,\alpha,E)$ is an asymptotic lift is a minor (and somewhat simpler) 
variation of the proof of 
the corresponding result of \cite{arvLift}.  Similarly, the proof of 
uniqueness is a straightforward variation of arguments in \cite{arvLift}; we omit the 
details. 

If a UP map $\phi:M\to M$ on a von Neumann algebra $M$ is completely positive, 
then its asymptotic lift $(N,\alpha,E)$ gives rise to 
a $W^*$-dynamical system $(N,\alpha)$ \cite{arvLift}.  While this need not be 
true for asymptotic lifts of UP maps, one can make the following assertion in general:

\begin{thm}\label{inthm2}
Let $\phi: M\to M$ be a UP map on a von Neumann algebra $M$ and 
let $(N,\alpha,E)$ be its asymptotic lift.  Then $N$ is order isomorphic 
to a unique JW$^*$-algebra in such a way that $\alpha$ is a Jordan automorphism of $N$.  
\end{thm}

Again, the proof follows along the lines of arguments in \cite{arvLift}, by 
introducing a Jordan multiplication on the range of a positive idempotent map as 
in Corollary 1.6 of \cite{EffSt}.  
The key step is to show that the for constructed asymptotic lift $(N,\alpha,E)$  
in which $N$ is the space of inverse sequences, there is a projection of 
norm one $Q: \ell^\infty(M)\to N$.  The existence 
of such a projection $Q$ follows from the argument of \cite{arvLift}.  Once one 
has a positive projection (which in the current setting 
is typically not {\em completely} positive) with 
these properties, 
one can introduce a Jordan project in $N$ by way of 
$$
x\circ y=Q(\frac{1}{2}(xy+yx)),\qquad x,y\in N, 
$$
and at this pont one can show that $N$ is order isomorphic to a JW$^*$-algebra 
by imitating arguments in \cite{arvLift}.  

\section{Embeddable asymptotic lifts}\label{S:em}

Let $(N,\alpha,E)$ be the asymptotic lift of a UP map $\phi:M\to M$ 
acting on a dual operator system.  
In this section we fix attention on those cases in which the asymptotic lift 
can be embedded as a subsystem of $M$ in the following particular way.  
We say that $(N,\alpha,E)$ {\em embeds in } $M$ 
if $E$ is restricts to an isometry on the self-adjoint part of $N$,
\begin{equation}\label{emEq1}
\|E(y)\|=\|y\|,\qquad y=y^*\in N.  
\end{equation}
This is equivalent to the assertion that $E$ implements 
an order isomorhism of $N$ onto $E(N)$.  Since all asymptotic lifts 
of $\phi:M\to M$ are isomorphic, this definition does not depend 
on the particular choice of $(N,\alpha,E)$.  
In such cases, the range $E(N)$ of $E$ is an operator subsystem 
of $M$ with the property that $\phi$ restricts to an 
automorphism of $E(N)$, 
$\alpha$ is identified 
with $\phi\restriction_{E(N)}$, and $E$ is 
identified with the inclusion map $\iota:E(N)\subseteq M$.  

The purpose of this section is to make a precise summary of those 
facts,  and especially to give an explicit description 
of $E(N)$ in terms of the action of $\phi$ on $M$.  
We begin by introducing the following 
operator system $M_\infty\subseteq M$: 
$$
M_\infty=\bigcup_{r>0}\bigcap_{n=1}^\infty \phi^n(\ball_r M).  
$$
$M_\infty$ consists of all $y\in M$ for which there is a 
{\em bounded} sequence $x_n\in M$ with $y=\phi^n(x_n)$, $n=1,2,\dots$.  
In general, $M_\infty$ is a self-adjoint linear subspace of $M$ 
containing the identity operator, it is invariant under $\phi$, 
and in fact  
$\phi(M_\infty)=M_\infty$.  Moreover, from (\ref{inEqB}) we may 
infer that 
\begin{equation}\label{emEq2}
E(N^{\rm{sa}})=\bigcup_{r>0}E(\ball_r N^{\rm{sa}})=
\bigcup_{r>0}\bigcap_{n=1}^\infty \phi^n(\ball_r M^{\rm{sa}})= M_\infty^{\rm{sa}},   
\end{equation}
hence $M_\infty=E(N)$ is precisely the range of $E$ in all cases.  We sometimes 
refer to $M_\infty$ as the {\em tail operator system} of $M$.  

\begin{rem}
It is clear that $M_\infty\subseteq \cap_n\phi^n(M)$, but 
the inclusion is typically proper.  As a simple example, let  
$M=M_2(\mathbb C)$ and consider the completely positive map $\phi:M\to M$ 
defined by 
$$
\phi
\begin{pmatrix}
a&b\\c&d
\end{pmatrix}
=
\begin{pmatrix}
a&\lambda b\\ \lambda c&d
\end{pmatrix}
$$
where $\lambda$ is a constant satisfying $0<\lambda<1$.  One has 
$\cap_n\phi^n(M)=M$, but in this case 
$M_\infty$ is the two-dimensional subalgebra 
of diagonal matrices.  
\end{rem}

\begin{prop}\label{emProp1}
Let $(N,\alpha,E)$ be an asymptotic lift of $\phi:M\to M$.  
Then $\ker E=\{0\}$ iff the restriction of $\phi$ to $M_\infty$ 
is both injective and surjective.  
\end{prop}

\begin{proof}    

In general, one has $\phi(M_\infty)=M_\infty$ by definition of $M_\infty$.  
Assuming that 
$\ker E=\{0\}$, choose $a\in M_\infty$ such that $\phi(a)=0$.  
Since $E(N)=M_\infty$, there is a $y\in N$ such that $a=E(y)$, hence 
$0=\phi(a)=\phi(E(y))=E(\alpha(y))$, and therefore $\alpha(y)=0$ because 
$E$ is injective.  $y=0$ follows because $\alpha$ is an automorphism of 
$N$, hence $a=E(y)=0$.  

Conversely, if the restriction of $\phi$ to $M_\infty$ is injective, 
choose a nonzero element $y\in N$.  Since the norms 
$\|E(\alpha^{-n}(y))\|$ increase to $\|y\|$ as $n\uparrow\infty$ (see Lemma 
3.7 of \cite{arvLift}), we must 
have $E(\alpha^{-n}(y))\neq 0$ for sufficiently large $n\geq 1$.  Since 
each power 
$\phi^n$ restricts to an injective map on $M_\infty=E(N)$, 
it follows that 
$$
E(y)=\phi^n(E(\alpha^{-n}(y)))\neq 0 
$$
for large $n$, hence $E$ is injective.   
\end{proof}

We conclude that whenever $\ker E=\{0\}$,  
the restriction of $\phi$ to $M_\infty$ defines 
an order-preserving linear bijection on $M_\infty$.  However, 
$M_\infty$ itself need 
not be closed in any topology in general, and the restriction of $\phi$ to $M_\infty$ need not 
be an order automorphism.  The following result implies that 
when $E$ restricts to an isometry on 
$N^{\rm{sa}}$, such anomalies cannot occur.  

\begin{thm}\label{emThm1}
For every UP map $\phi: M\to M$ on a dual operator system, 
the following are equivalent.
\begin{enumerate}
\item[(i)] The asymptotic lift of $\phi:M\to M$ embeds in $M$.  
\item[(ii)] $M_\infty$ is weak$^*$-closed and 
$\phi$ restricts to an order automorphism of it.  
\item[(iii)] The asymptotic lift of $\phi:M\to M$ is isomorphic to 
the triple $(M_\infty, \phi\restriction_{M_\infty},\iota)$, where 
$\iota: M_\infty\subseteq M$ is the inclusion map.  
\end{enumerate}
\end{thm}

\begin{proof}
Let $(N,\alpha,E)$ be an asymptotic lift of $\phi:M\to M$.  

(i)$\implies$(ii):  By hypothesis, $E$ restricts to an 
isometry of $N^{\rm {sa}}$ onto the self-adjoint part of 
$
E(N)=M_\infty
$.  
Since $M_\infty$ is a self-adjoint linear space of operators and 
the adjoint operation is weak$^*$-continuous, 
$M_\infty$ will be weak$^*$-closed provided we show that its self-adjoint 
part $M_\infty^{\rm{sa}}$ is weak$^*$-closed.  By a standard result on 
the weak$^*$-closure of convex sets in dual Banach spaces, 
this will follow if 
we show that for every $r>0$, the intersection $M_\infty^{\rm{sa}}\cap\ball_r M$ 
of $M_\infty^{\rm{sa}}$ with the $r$-ball of $M$ 
is weak$^*$-closed.  
Now the connecting map $E: N\to M$ restricts to a 
weak$^*$-continuous isometry of $N^{\rm{sa}}$ onto $M_\infty^{\rm{sa}}$, 
and by (\ref{inEqB}),  it carries $\ball_r N^{\rm{sa}}$ onto 
$M_\infty^{\rm{sa}}\cap\ball_r M$.  Since the 
unit ball of $N$ is weak$^*$-compact, it follows that 
$M_\infty^{\rm{sa}}\cap\ball_r M=E(\ball_r N^{\rm{sa}})$ is weak$^*$-compact, hence 
weak$^*$-closed.  We conclude that $M_\infty^{\rm{sa}}$ is 
weak$^*$-closed.    

Making use of Proposition \ref{emProp1}, let $\psi$ be the linear automorphism of $M_\infty$ inverse to 
$\phi\restriction_{M_\infty}$.  Then since $E\circ\alpha=\phi\circ E$, 
we have 
$E\circ\alpha^{-1}=\psi\circ E$.  Since both $E$ and $\alpha^{-1}$ restrict to 
isometries on $N^{\rm{sa}}$, it follows that $\psi$ restricts to an 
isometry on $M_\infty^{\rm{sa}}=E(N^{\rm{sa}})$.  Since $\psi(\mathbf 1)=\mathbf 1$, 
$\psi$ must also be order-preserving, and we conclude that $\phi\restriction_{M_\infty}$ is 
an order automorphism.  

(ii)$\implies$(iii): Under the hypothesis (ii), $(M_\infty,\phi\restriction_{M_\infty}, \iota)$ 
becomes a reversible lift of $\phi:M\to M$ that is obviously nondegenerate. 
We show that it is the asymptotic lift of $\phi:M\to M$ by establishing (\ref{inEqA}).  
For that, consider the decreasing sequence of weak$^*$-compact sets 
$$
\phi^n(\ball M^{sa}),\qquad n=1,2,\dots  
$$
and choose a normal linear functional $\rho\in M_*$.  By (\ref{inEqB}), we have 
$$
\bigcap_{n=1}^\infty \phi^n(\ball M^{sa})=E(\ball N^{\rm{sa}}).    
$$
Since $\rho$ is weak$^*$-continuous, it follows from Lemma 3.5 of \cite{arvLift}) that 
\begin{align*}
\|\rho\circ \phi^n\restriction_{M_\infty^{\rm{sa}}}\|&=
\sup\{\rho(x): x\in \phi^n(\ball M^{\rm{sa}})\}\downarrow 
\sup\{\rho(x): x\in E(\ball N^{\rm{sa}})\}
\\
&=\|\rho\circ E\restriction_{N^{\rm{sa}}}\|,   
\end{align*}
as $n\uparrow \infty$, and (iii) follows.  

The implication (iii)$\implies$(i) is obvious.   
\end{proof}

Taken together, theorems \ref{emThm1} and \ref{inthm2} imply: 

\begin{cor}\label{emCor1}
Let $\phi: M\to M$ be a UP map on a von Neumann algebra $M$ 
whose asymptotic lift embeds in $M$.  
Then the tail operator system 
$M_\infty$ is order-isomorphic to a JW$^*$-algebra in such a way 
that the restriction of $\phi$ to $M_\infty$ becomes a Jordan automorphism, 
and the asymptotic lift of $\phi$ is the triple 
$(M_\infty,\phi\restriction{M_\infty},\iota)$, $\iota$ being the inclusion map 
$\iota:M_\infty\subseteq M$.  
\end{cor}

The asymptotic lifts of many UP maps are embeddable.  That is true for 
all examples covered by the hypotheses of \cite{arvStableI}, 
and in Section \ref{S:tm} below we elaborate on the special case of maps 
on finite-dimensional algebras.  
Not all infinite-dimensional UP maps are embeddable, and 
the following example illustrates the fact.  

\begin{rem}[An Example]
Let $L^2=L^2(\mathbb T,\frac{d\theta}{2\pi})$, let $H^2\subseteq L^2$ 
be the usual Hardy space, and let $s\in\mathcal B(H^2)$ be the unilateral shift.  
Consider the UCP map $\phi$ defined on $M=\mathcal B(H^2)$ by 
$$
\phi(x)=s^*xs,\qquad A\in \mathcal B(H^2).  
$$
Let $u\in\mathcal B(L^2)$ be the bilateral shift and let $p$ be the 
projection of $L^2$ on $H^2$.  Then the asymptotic lift of $\phi$ is 
the triple $(\mathcal B(L^2),\alpha, E)$, where $\alpha(y)=u^*bu$ and 
$E$ is the compression map $E(y)=py\restriction_{H^2}$.  Since $\ker E\neq\{0\}$, 
the asymptotic lift of $\phi$ is not embeddable.  
We omit the calculations.  
\end{rem}

\section{Embedding and the multiplicative core}\label{S:ec}

If the asymptotic lift of a UP map $\phi: M\to M$ acting on 
a von Neumann algebra embeds in $M$, then by Corollary \ref{emCor1}, 
$M_\infty$ is order isomorphic to a JW$^*$-algebra and the restriction 
of $\phi$ to $M_\infty$ becomes a Jordan automorphism of $M_\infty$.  
Notice that this does not imply that $M_\infty$ itself is closed 
under the Jordan product  $x\circ y=\frac{1}{2}(xy+yx)$ inherited 
from $M$.    
Section \ref{S:tm} contains a simple example that illustrates 
the phenomenon.  

In the present section we study the Jordan 
structure of $M_\infty$ in more detail in the case when 
the asymptotic lift embeds in $M$.  In particular, we show that if 
$M_\infty$ is itself a JW$^*$-algebra then it coincides with 
the {\em multiplicative core} $C_\phi$ of $\phi$ 
introduced in \cite{storPos}.  We also give necessary and sufficient conditions 
for $(C_\phi,\phi\restriction_{C_\phi},\iota)$ to be the asymptotic lift.

Following  \cite{storPos}, the 
{\em definite set} of a UP map $\phi:M\to M$ is defined by 
\begin{equation}\label{mcEq1}
M_\phi=\{x\in M: \phi(x^*\circ x)=\phi(x)^*\circ \phi(x)\};
\end{equation}
it is a JW$^*$-algebra, and we have 
$$
\phi(x\circ y)=\phi(x)\circ \phi(y), \qquad x\in M_\phi,\quad y\in M.     
$$
Continuing as in \cite{storPos}, 
one can show that  
$$
B_\phi=\{x\in M:  \phi^n(x)\in M_\phi,\quad n=0,1,2,\dots\}   
$$
is a JW$^*$-algebra on which $\phi$ restricts to a Jordan endomorphism, 
and one forces surjectivity on $\phi$ by restricting it to the 
``tail" subalgebra 
$$
C_\phi=\bigcap_{n=1}^\infty\phi^n(B_\phi).  
$$
$C_\phi$ is called the {\em multiplicative core} of $\phi$.  
We require the following variation of Lemma 6 of \cite{storPos}.  

\begin{lem}\label{ecLem1}
The multiplicative core is characterized as 
the largest JW$^*$-subalgebra $N\subseteq M$ with the following two 
properties:
\begin{enumerate}
\item[(i)] $\phi(N)= N$.  
\item[(ii)]$\phi(x\circ y)=\phi(x)\circ \phi(y)$ for all $x,y\in N$.  
\end{enumerate}
\end{lem}

\begin{proof}
Let $N$ be a JW$^*$-subalgebra of $M$ with properties (i) and (ii), 
and choose $x\in N$.  Then we have 
$$
\phi(\phi^n(x)^*\circ\phi^n(x))=\phi^{n+1}(x)^*\circ\phi^{n+1}(x),\qquad n\geq 0,
$$
because $\phi$ is a Jordan endomorphism of $N$, hence $N\subseteq B_\phi$.  
It follows that 
$N=\phi^n(N)\subseteq \phi^n(B_\phi)$ for every $n\geq 0$, 
so that $N\subseteq \cap_n\phi^n(B_\phi)=C_\phi$.  
\end{proof}

The 
following result associates the multiplicative core 
with asymptotic lifts that 
are embeddable as Jordan subalgebras:

\begin{thm}\label{ecThm1}
For every UP map $\phi:M\to M$ acting on a von Neumann algebra 
$M$,  the tail operator system $M_\infty$ contains the multiplicative core $C_\phi$.  

Assume that the asymptotic lift embeds in $M$.  Then 
$M_\infty$ is a Jordan subalgebra of $M$ precisely when $M_\infty=C_\phi$.  
\end{thm}

\begin{proof}
By definition, $C_\phi=\cap_n\phi^n(B_\phi)$, and 
$\phi$ restricts to a Jordan endomorphism of $B_\phi$.  Since 
each power $\phi^n$ of $\phi$ restricts to  a Jordan 
endomorphism on $B_\phi$ and a Jordan homomorphism maps the unit ball of the 
self-adjoint part of its domain onto the unit ball of the self-adjoint part of 
its range, we have 
$\phi^n(\ball_r B_\phi^{\rm{sa}})=\ball_r(\phi^n(B_\phi)^{\rm{sa}})$ 
for every $r>0$ and every $n=1,2,\dots$, hence 
$$
\ball_r C_\phi^{\rm{sa}}=\bigcap_{n=1}^\infty\ball_r\phi^n(B_\phi^{\rm{sa}}) 
=\bigcap_{n=1}^\infty\phi^n(\ball_r B_\phi^{\rm{sa}})\subseteq M_\infty.   
$$ 
The asserted inclusion 
$C_\phi=C_\phi^{\rm{sa}}+iC_\phi^{\rm{sa}}\subseteq M_\infty$ 
follows after taking the union over all positive $r$.  

Assume that the asymptotic lift of $\phi:M\to M$ embeds in $M$.  
If  $M_\infty$ 
is closed under the Jordan multiplication of $M$, then by Theorem \ref{emThm1} (iii), 
$\phi\restriction_{M_\infty}$ is 
an order automorphism of $M_\infty$, and an application of Kadison's Schwarz 
inequality for positive linear maps implies that $\phi$ is a Jordan 
automorphism of $M_\infty$.  By Lemma \ref{ecLem1}, 
$M_\infty\subseteq C_\phi$, and  
we conclude that $M_\infty=C_\phi$.  The converse is trivial.  
\end{proof}

It is significant that when 
$M_\infty$ is not itself closed under the Jordan product, it contains 
a largest JW$^*$-algebra that is characterized as follows.  

\begin{prop}\label{mcProp1}
Let $\phi:M\to M$ be a UP map whose asymptotic lift embeds in $M$.  Then the weak$^*$-closed 
linear span of all projections in $M_\infty$ is a JW$^*$-algebra.  
\end{prop}

\begin{proof}
By Corollary \ref{emCor1}, there is a JW$^*$-algebra $B$ and an order isomorphism 
$\alpha$ of $B$ onto $M_\infty$.  Viewing $\alpha$ as a UP map of $B$ into the 
JW$^*$-algebra $M$, it makes sense to speak of its definite set 
$B_\alpha$.  Since $B_\alpha$ is a JW$^*$-subalgebra of $B$ and the restriction 
of $\alpha$ to $B_\alpha$ is  a Jordan homomorphism, 
it follows that $\alpha(B_\alpha)\subseteq M_\infty$ is a JW$^*$-subalgebra of $M$.  

Being a JW$^*$-algebra, $\alpha(B_\alpha)$ is generated by its projections.  
Conversely, we claim that every 
projection  $e\in M_\infty$ belongs to $\alpha(B_\alpha)$.  To see that, fix 
such an $e$ and 
choose $f\in B$ such that $\alpha(f)=e$.  Then $0\leq f\leq \mathbf 1$ because 
$\alpha$ is a UP order isomorphism, so by Kadison's Schwarz inequality 
$$
0=e-e^2=\alpha(f)-\alpha(f)^2\geq \alpha(f-f^2)\geq 0,   
$$
hence $\alpha(f-f^2)=0$ and finally $f=f^2$.  Since $\alpha(f^2)=\alpha(f)^2$, 
$f$ belongs to the definite set of $\alpha$, hence $e=\alpha(f)\in\alpha(B_\alpha)$.  
\end{proof}

We conclude the section by describing intrinsic conditions on a UP map which imply that its 
asymptotic lift is embeddable as the multiplicative core; note that the sufficient conditions 
of Theorem \ref{mcThm2} are clearly necessary as well.

For a linear functional $\rho$ on $M$ and a set of operators $S\subseteq M$, 
we write $\rho\perp S$ whenever $\rho(S)=\{0\}$.  We require the following result characterizing 
the equality $M_\infty=C_\phi$ in terms of the action of $\phi$ on the predual of $M$:

\begin{lem}\label{mcLem2}
Let $\phi:M\to M$ be a UP map.  Then $M_\infty=C_\phi$ iff 
for every $\rho\in M_*$ satisfying $\rho\perp C_\phi$, one has  
\begin{equation}\label{mcEq2}
\lim_{n\to\infty}\|\rho\circ \phi^n\|=0.  
\end{equation}
\end{lem}

\begin{proof}
We claim first that $\rho\in M_*$ satisfies (\ref{mcEq2}) iff $\rho\perp M_\infty$.  Indeed, 
setting $M_\infty(r)=\cap_n\phi^n(\ball_r M)$ for $r>0$, and noting that the 
compact convex sets $\phi^n(\ball_r M)$ decrease to $M_\infty(r)$ as $n\uparrow\infty$, 
we can 
apply Lemma 3.5 of \cite{arvLift} to conclude that (\ref{mcEq2}) is equivalent to the assertion 
\begin{align*}
\sup\{|\rho(x)|: x\in M_\infty(r)\}&=\lim_{n\to\infty}\sup\{|\rho(x)|: x\in \phi^n(\ball_r M)\}
\\
&=r\cdot \lim_{n\to\infty}\|\rho\circ\phi^n\|=0, 
\end{align*}
for every $r>0$.  Noting that $M_\infty=\cup_{r>0}M_\infty(r)$, the claim follows.

The preceding paragraph, together with a standard separation theorem, shows that the 
assertion 
\begin{equation}\label{mcEq3}
\rho\perp C_\phi\implies \lim_{n\to\infty}\|\rho\circ\phi^n\|=0
\end{equation}
is equivalent 
to the assertion $\overline{M_\infty}^{{\rm w}*}\subseteq C_\phi$; and since in 
general we have $C_\phi\subseteq M_\infty\subseteq \overline{M_\infty}^{{\rm w}*}$, 
(\ref{mcEq3}) is seen to be equivalent to $C_\phi=M_\infty$.  
\end{proof}

\begin{thm}\label{mcThm2}
Let $\phi:M\to M$ be a UP map such that 
\begin{enumerate}
\item[(i)]
The positive linear map obtained by restricting $\phi$ to $C_\phi$ is faithful.  
\item[(ii)] $\lim_{n\to\infty}\|\rho\circ\phi^n\|=0$ for every $\rho\in M_*$ satisfying 
$\rho\perp C_\phi$.  
\end{enumerate}
Then $M_\infty=C_\phi$, the restriction of 
$\phi$ to $C_\phi$ is a Jordan automorphism of $C_\phi$, 
and the asymptotic lift 
of $\phi$ is $(C_\phi,\phi\restriction_{C_\phi},\iota)$.   
\end{thm}

\begin{proof}
Hypothesis (i) implies that the restriction of $\phi$ to 
$C_\phi$ is a Jordan automorphism.  Thus, the triple $(C_\phi,\phi\restriction_{C_\phi},\iota)$ 
is a nondegenerate reversible lift of $\phi$.  To show that it is the asymptotic lift, 
we must establish the following inequality for every $\rho\in M_*$:
\begin{equation}\label{mcEq4}
\lim_{n\to\infty}\|\rho\circ\phi^n\|\leq \|\rho\restriction_{C_\phi}\|.  
\end{equation}
For that, fix $\rho$ and, for every $n=1,2,\dots$, choose an element $x_n\in M$ 
satisfying $\|x_n\|= 1$ and $|\rho(\phi^n(x_n))|=\|\rho\circ\phi^n\|$.  We can 
find a subsequence $n_1<n_2<\dots$ such that $\phi^{n_k}(x_{n_k})$ converges to $y\in M$ 
weak$^*$ as $k\to\infty$.  
Since the weak$^*$-compact sets $\phi^n(\ball M)$ decrease with increasing $n$, $y$ must belong 
to their intersection 
$\cap_n\phi^n(\ball M)=M_\infty$, and of course $\|y\|\leq 1$.  

Making use of hypothesis 
(ii) and Lemma \ref{mcLem2}, we conclude that $M_\infty=C_\phi$, so that $y\in \ball C_\phi$.  
Hence 
$$
\lim_{n\to\infty}\|\rho\circ\phi^n\|=
\lim_{k\to\infty}\|\rho\circ\phi^{n_k}\|=\lim_{k\to\infty}|\rho(\phi^{n_k}(x_{n_k}))|=|\rho(y)|
\leq\|\rho\restriction_{C_\phi}\|, 
$$
and the desired inequality (\ref{mcEq4}) follows.  
\end{proof}

\begin{rem}
The conclusion of Theorem \ref{mcThm2} can be significantly strengthened whenever 
there is a normal positive projection $E$ of $M$ onto $C_\phi$ - in 
particular, whenever $M$ is a finite von Neumann algebra.  In such cases, 
a simple argument (that we omit) shows that the automorphism 
$\alpha=\phi\restriction_{C_\phi}$ of $C_\phi$ satisfies 
$$
\lim_{n\to\infty}\|\rho\circ\phi^n-\rho\circ\alpha^n\circ E\|=0,\qquad \rho\in M_*.   
$$
\end{rem}

\section{UP maps on finite-dimensional algebras}\label{S:tm}

In this section we show that the asymptotic lift of every UP map acting on a 
finite-dimensional algebra $M$ embeds in $M$.  We identify 
$M_\infty$ with the multiplicative core whenever there is a faithful 
$\phi$-invariant state, and more generally, we identify the multiplicative core 
when $\phi$ is faithful.  We conclude with an elementary example exhibiting nontrivial 
asymptotic dynamics, for which 
$M_\infty$ is not closed under the Jordan product of $M$ and hence differs from the 
multiplicative core.

\begin{thm}\label{tmThm1}
Let $\phi:M\to M$ be a UP map on a finite-dimensional von Neumann algebra.  
Then $(M_\infty,\phi\restriction_{M_\infty}, \iota)$ is the 
asymptotic lift of $\phi$.  If, in addition, for every positive 
operator $x\in M$ one has 
\begin{equation}\label{tmEq1}
\lim_{n\to\infty}\|\phi^n(x)\|=0\implies x=0, 
\end{equation}
then $M_\infty$ is the multiplicative core $C_\phi$.  Condition 
(\ref{tmEq1}) is satisfied whenever there is a faithful state $\rho$ of $M$ 
satisfying $\rho\circ\phi=\rho$.   
\end{thm}

The proof requires a known elementary result (see \cite{EffSt}): 

\begin{lem}\label{tmLem1}
Let $M$ be a unital JW$^*$ algebra  and 
let $E:M\to M$ be an idempotent UP map that is faithful:
$x\in M^+$, $E(x)=0\implies x=0$.   
Then $E(M)$ is a Jordan subalgebra of $M$.
\end{lem}

\begin{proof}
Choose a self-adjoint element $x\in E(M)$.  
Then $E(E(x^2)-x^2)=E(x^2)-E(x^2)=0$. By Kadison's Schwarz inequality, 
$x^2=E(x)^2\leq E(x^2)$, so that $E(x^2)-x^2\geq 0$.   Since $E$ is faithful, 
$E(x^2)=x^2\in E(M)$.  
This shows that $E(M)^{\rm{sa}}$ is closed under the Jordan product, hence $E(M)$ 
is a Jordan subalgebra of $M$.  
\end{proof}

\begin{proof}[Proof of Theorem \ref{tmThm1}]
There is a sequence $n_1<n_2<\cdots$ of positive integers 
such that $\phi^{n_k}$ converges to a unique idempotent $E$ (this is a result 
of Kuperberg, see 
Theorem 4.1 of \cite{arvStableI} et. seq.).  Note that 
$E(M)=M_\infty$.  Indeed, for every $x\in M$, 
$$
E(x)=\lim_k\phi^{n_k}(x)\in\cap_n\overline{\{\phi^n(x),\phi^{n+1}(x)\dots\}}\subseteq \cap_n\phi^n(\ball_{\|x\|} M),  
$$  
hence 
$E(x)\subseteq M_\infty$.  For the opposite inclusion, 
choose $y\in M_\infty$.  Then there is a bounded sequence $x_k\in M$ such 
that $y=\phi^{n_k}(x_k)$ for every $k$.  Let $k^\prime$ be a subseqence of $k$ such 
that $x_{{k^\prime}}$ converges to $x\in M$.  Then 
$$
\|y-\phi^{n_{k^\prime}}(x)\|\leq \|y-\phi^{n_{k^\prime}}(x_{k^\prime})\|+
\|\phi^{n_{k^\prime}}(x_{k^\prime})-\phi^{n_{k^\prime}}(x)\|\leq \|x_{k^\prime}-x\|, 
$$
and the right side tends to zero as $k^\prime\to \infty$.  

It was also shown in \cite{arvStableI} that $\phi$ restricts to a surjective isometry 
on $M_\infty$ and that the powers of $\phi$ tend to 
zero on $\ker E$, so that 
$$
\lim_{n\to\infty}\|\phi^n\circ E-\phi^n\|=0.  
$$
It follows that for every bounded linear functional $\rho$ on $M$, we have 
\begin{align*}
\limsup_{n\to\infty}|\,\|\rho\restriction_{M_\infty}\|-\|\rho\circ\phi^n\|\, |
&=\limsup_{n\to\infty}|\,\|\rho\circ\phi^n\circ E\|-\|\rho\circ\phi^n\|\, |
\\
&\leq\lim_{n\to\infty}\|\rho\circ\phi^n\circ E-\rho\circ\phi^n\|=0.    
\end{align*}
The latter implies that $(M_\infty,\phi\restriction_{M_\infty},\iota)$ 
is the asymptotic lift of $\phi$.  

Assuming that (\ref{tmEq1}) is satisfied, we claim that $M_\infty=C_\phi$.  
By Theorem \ref{ecThm1},  
it is enough to show that $M_\infty$ is closed under the Jordan multiplication 
of $M$.  
Note that $E$ is must be faithful.  
Indeed, if $x\geq 0$ and $E(x)=0$, then 
$$
\lim_k\|\phi^{n_k}(x)\|=\|E(x)\|=0, 
$$
and since the sequence of norms $\|\phi_n(x)\|$ is decreasing, (\ref{tmEq1}) 
implies that $x=0$.  By Lemma \ref{tmLem1}, $M_\infty$ is a Jordan subalgebra 
of $M$.  

Finally, assuming that there is a faithful 
$\phi$-invariant state $\rho$, we claim that (\ref{tmEq1}) holds.  
Indeed, if $x$ is a positive operator satisfying 
$\|\phi^n(x)\|\to 0$ as $n\to \infty$, then 
$$
|\rho(x)|=|\rho\circ\phi^n(x))|\leq \|\phi^n(x)\|\to 0
$$
as $n\to\infty$.  Hence $\rho(x)=0$, 
and $x=0$ follows because $\rho$ is faithful.  
\end{proof}

We conclude with the following result, which 
identifies the multiplicative core $C_\phi$ in many cases where $C_\phi\neq M_\infty$.  

\begin{cor}
Let $\phi:M\to M$ be a faithful UP map on a finite-dimensional von Neumann algebra.  
Then we have:
\item [(i)]  The multiplicative core $C_{\phi}$ of $\phi$ is 
the linear space spanned by the projections in $M_{\infty}$.  
\item [(ii)]   $M_{\infty}=C_{\phi}$ if and only if there exists a $\phi$ -invariant
 state of $M$ whose restriction to the JW$^*$ algebra generated by $M_{\infty}$ is faithful.
\end{cor}

\begin{proof} (i). Let $A$ be the linear span of all projections in $M_\infty$.  
We claim that $\phi(A)=A$.  Indeed, if $e$ is a projection in $A$ 
then since $\phi$ is an order automorphism of $M_\infty$ 
there is an operator $f\in M_\infty$ such that $0\leq f\leq\mathbf 1$ 
and $\phi(f)=e$.  As in the proof of Proposition \ref{mcProp1}, this implies 
$\phi(f-f^2)=0$, hence $f=f^2$ is a projection because $\phi$ is faithful.  
This implies that $A\subseteq \phi(A)$ and, since $A$ is finite-dimensional, 
$\phi(A)=A$.  

Proposition \ref{mcProp1} implies that $A$ is the largest Jordan algebra 
in $M_\infty$.  
The multiplicative core $C_\phi$ is a Jordan 
algebra in $M_\infty$ by Theorem \ref{ecThm1}, hence $C_\phi\subseteq A$.    Lemma 
\ref{ecLem1} now implies that $C_\phi=A$. 

(ii) Let $E$ be the $\phi$-invariant projection of $M$ onto $M_{\infty}$ described in the 
proof of Theorem (\ref{tmThm1}), and let $P $ denote the $\phi$-invariant
 projection of $C_{\phi}$ onto the set $M^{\phi}$ of fixed points of $\phi$, (see \cite{EffSt}).
 Since $\phi$ is a Jordan automorphism of $C_{\phi}$ $P$ is faithful, so that $M^{\phi}$ is a JW{$^*$}
 algebra by Lemma (\ref{tmLem1}). Suppose $\omega$ is a faithful state of $M^{\phi}$.  Then if 
 $M_{\phi}=C_{{\phi}}$ the state $\rho=\omega \circ P\circ E$ is the desired state on $M$. 
 
 Conversely if $\rho$ is a $\phi$-invariant state whose restriction to the JW{$^*$}algebra $N$
  generated by $M_{\infty}$ is faithful, then $\rho= \rho\circ E$, so the restriction of $E$
  to $N$ is faithful.  Thus by Lemma (\ref{tmLem1}) $M_{\infty}$ is a JW{$^*$} algebra, hence by 
  Theorem (\ref {ecThm1}) equal to $C_{\phi}$.
  \end{proof}

We conclude by describing an example of a UP map 
on the $3$-dimensional commuative 
\cstar\ $M=\mathbb C^3$ for which $M_\infty$ is not closed under the ambient 
Jordan multiplication of $M$.  While there are simpler examples with 
that specific property, this one 
exhibits nontrivial asymptotic dynamics that are 
not detected by the multiplicative core.  
Viewing the elements of $M$ as column vectors, 
the map $\phi$ is multiplication by the stochastic matrix
$$
\phi=
\begin{pmatrix}
\frac{1}{3}&\frac{1}{3}&\frac{1}{3}\\
0&0&1\\
0&1&0
\end{pmatrix}
.
$$
The even and odd powers of $\phi$ are 
$$
\phi^{2n}=
\begin{pmatrix}
\frac{1}{9^n}&\frac{1}{2}-\frac{1}{2\cdot 9^n}&\frac{1}{2}-\frac{1}{2\cdot 9^n}\\
0&1&0\\
0&0&1
\end{pmatrix}
, 
\  
\phi^{2n+1}=
\begin{pmatrix}
\frac{1}{3\cdot 9^n}&\frac{1}{2}-\frac{1}{6\cdot 9^n}&\frac{1}{2}-\frac{1}{6\cdot 9^n}\\
0&0&1\\
0&1&0
\end{pmatrix}
,
$$
and the unique {\em idempotent} limit point of $\{\phi,\phi^2,\phi^3,\dots\}$ is 
given by 
$$
E=\lim_{n\to\infty}\phi^{2n}=
\begin{pmatrix}
0&\frac{1}{2}&\frac{1}{2}\\
0&1&0\\
0&0&1
\end{pmatrix}
.
$$
The range of $E$ is the two-dimensional space (written as  
row vectors) 
$$
M_\infty=E(M)=\{(\frac{a+b}{2},a,b): a,b\in\mathbb C\}.
$$ 

%\begin{prop}\label{tmProp1}
We summarize the basic properties of this example without proof:

\noindent
{\em 
Relative to the intrinsic (Jordan) multiplication defined by 
$x\circ y=E(xy)$, 
$M_\infty$ is isomorphic to the two-dimensional commutative 
\cstar\ $\mathbb C^2$.  This identification 
implements a conjugacy of 
$\phi\restriction_{M_\infty}$ with the order $2$ automorphism 
$(a,b)\mapsto (b,a)$, $a,b\in\mathbb C$.  
The multiplicative core of $\phi$ is the one-dimensional \cstar\  $C_\phi=\mathbb C\cdot 1$.  
}
%\end{prop}

\begin{rem}[Invariant states of $\phi$]
Perhaps it is worth pointing out 
that the above UP map $\phi:\mathbb C^3\to\mathbb C^3$ 
has a unique invariant state, namely  
$$
\rho(a,b,c)=\frac{b+c}{2},\qquad (a,b,c)\in\mathbb C^3.  
$$
This state is not faithful of course, so there is no conflict 
with Theorem \ref{tmThm1}.  
\end{rem}

%\vfill
%\eject

%\bibliography{bibData}  %Remove this line when finished, see below.  
\bibliographystyle{alpha}

\newcommand{\noopsort}[1]{} \newcommand{\printfirst}[2]{#1}
  \newcommand{\singleletter}[1]{#1} \newcommand{\switchargs}[2]{#2#1}

\end{document}